\documentclass[12pt]{article}

\usepackage[dvips]{graphicx,psfrag}
\usepackage{floatflt,epsfig,epsf}

\begin{document}

\newcommand{\pf}{

\smallskip

\noindent {\it Proof : }}

\newcommand{\nul}{\mbox{$\tilde{0}$}}

\newcommand{\al}{\mbox{$\alpha$}}

\newcommand{\be}{\mbox{$\beta$}}

\newcommand{\dd}{\mbox{$\delta$}}
\newcommand{\la}{\mbox{$\lambda$}}

\newcommand{\e}{\mbox{$\epsilon$}}

\newcommand{\norm}[1]{\mbox{$\|#1\|$}}

\newcommand{\cvf}{\mbox{$\stackrel{w}{\rightarrow}$}}

\newcommand{\hs} {\hskip}

\newcommand{\vs} {\vskip}

\newcommand {\ram}{[\omega]^\omega}

\newcommand {\ca} {2^\omega}

\newcommand{\cantor}{2^{<\omega}}

\newcommand {\N}{\mathbb N}
\newcommand {\poset}{{\mathcal P}}
\newcommand {\Q}{\mathbb Q}

\newcommand {\R}{\mathbb R}

\newcommand {\Z}{\mathbb Z}

\newcommand {\C}{\mathbb C}

\newcommand{\om}{\omega}

\newcommand{\eps}{\epsilon}

\newcommand{\compl}{\complement}

\newcommand {\iso} {\cong}

\newcommand{\isom}{\simeq}

\newcommand{\ext}{\preccurlyeq}

\newcommand {\tom} {\emptyset}

\newcommand{\emb}{\sqsubseteq}

\newcommand{\inj}{\hookrightarrow}

\newcommand{\begr}{\upharpoonright}

\newcommand{\surj}{\twoheadrightarrow}

\newcommand{\bij}{\longleftrightarrow}

\newcommand{\hviss}{\longleftrightarrow}

\newcommand{\hvis}{\Longleftarrow}

\newcommand{\saa}{\Longrightarrow}

\newcommand{\til}{\longrightarrow}

\newcommand{\equi}{\Longleftrightarrow}

\newcommand{\Lim}[1]{\mathop{\longrightarrow}\limits_{#1}}

\newcommand {\Del}{ \; \Big| \;}

\newcommand {\del}{ \; \big| \;}

\newcommand {\Mgdv}{\Big\{}

\newcommand {\mgdv}{\big\{}

\newcommand {\Mgdh}{\Big\}}

\newcommand {\mgdh}{\big\}}

\newcommand {\Intv}{\Big[}

\newcommand {\intv}{\big[}

\newcommand {\Inth}{\Big]}

\newcommand {\inth}{\big]}

\newcommand {\For}{\Bigcup}

\newcommand {\for}{\bigcup}

\newcommand {\Snit}{\Bigcap}

\newcommand {\snit}{\bigcap}

\newcommand {\og}{\; \land \;}

\newcommand {\eller}{\; \vee\;}

\newcommand{\ikke}{\lnot}

\newcommand {\go} {\mathfrak}

\newcommand {\ku} {\mathcal}

\newcommand {\un} {\underline}

\newcommand{\hr} {\hookrightarrow}

\newcommand{\st} {\stackrel}

\newcommand{\cqd} {\hspace{10pt} \rule {5pt}{5pt}}

\newcommand {\ex} {\exists}

\renewcommand {\a} {\forall}

\newcommand{\fed}{\mathbf}

\font\tenBbb=msbm10  \font\sevenBbb=msbm7  \font\fiveBbb=msbm5

\newfam\Bbbfam

\textfont\Bbbfam=\tenBbb \scriptfont\Bbbfam=\sevenBbb

\scriptscriptfont\Bbbfam=\fiveBbb

\def\Bbb{\fam\Bbbfam\tenBbb}

\def\R{{\Bbb R}}

\def\C{{\Bbb C}}

\def\N{{\Bbb N}}

\def\zed{{\Bbb Z}}

\def\Q{{\Bbb Q}}

\newcommand{\pff}{$ $}

\newtheorem{defi}{Definition}

\newtheorem{prop}[defi]{Proposition}

\newtheorem{lemm}[defi]{Lemma}

\newtheorem{theo}[defi]{Theorem}

\newtheorem{coro}[defi]{Corollary}

\newtheorem{rema}[defi]{Remark}

\newtheorem{conj}[defi]{Conjecture}

\title{ Some equivalence relations which are Borel reducible to 
 isomorphism  between separable
Banach spaces}
\date{ }
\author{Valentin Ferenczi\footnote{This author was supported by FAPESP
 Grant 2002/09662-1.}\ \ and El\'oi Medina Galego}
\maketitle

\begin{abstract}{\em We prove that the relation $E_{K_{\sigma}}$ is
 Borel
    reducible to isomorphism and
 complemented biembeddability between subspaces of $c_0$ or $l_p$ with
 $1 \leq
 p <2$. We also show that the relation
$E_{K_{\sigma}} \otimes =^+$ is Borel reducible to isomorphism,
 complemented biembeddability, and Lipschitz isomorphism between
 subspaces of
 $L_p$ for $1 \leq p <2$.}
\end{abstract}

\vs 0.2cm

{\bf 1. Introduction}

\vs 0.2cm

In this paper, we are mainly interested in the complexity of the
 relation of
isomorphism between separable Banach spaces.  The central notion in
 the theory
of classification of analytic equivalence relations on Polish spaces
 by means
of their relative complexity is the concept of Borel reducibility
 between
equivalence relations. This concept originated from the works of H.
 Friedman and L. Stanley and independently
 from the works of L. A. Harrington, A. S. Kechris and A. Louveau.\footnote {2000 {\it {Mathematics
 Subject Classification.}} Primary 03E15, 46B03. {\it {Key words and
 phrases:}} Analytic and Borel equivalence relations, Complexity of
 the relation of isomorphism between separable Banach spaces}

\vs 0.2cm

{\bf 1.1. Borel reducibility of equivalence relations on Polish
 spaces.} Let $R$ (resp. $R'$) be an analytic equivalence relation on
 a Polish space $E$
(resp. $E'$). We say that $(E,R)$ is {\em Borel reducible} to
 $(E',R')$, and write $(E,R) \leq_B (E',R')$, if there exists a Borel
 map $f$ from $E$ to $E'$,
such that for all $x$ and $y$ in $E$,
$$ x R y \Leftrightarrow f(x) R' f(y).$$

One may also restricts one's attention to Borel, instead of analytic,
equivalence relations. Note also that the above definition is valid
 for
general relations. A theory of $\leq_B$ for quasi-orders has been
 recently developed by A. Louveau and C. Rosendal \cite{LR}.

Observe that if $R$ and $R'$ are equivalence relations, then the
 induced quotient
map from $E/R$ into $E'/R'$ is an injection. In particular, $E'$ has
 at least as many
$R'$-classes as $E$ has $R$-classes. In fact, equivalence classes for
 $R'$ provide invariants for the equivalence relation $R$, and
 furthermore this can be obtained in a Borel way. So the order
 $\leq_B$ can be seen as  a  measure of relative complexity between
 analytic equivalence relations on Polish spaces. 

The relation $(E,R)$ is {\em Borel bireducible} to $(E',R')$,
 $(E,R) \sim_B (E',R')$,  whenever 
both $(E,R) \leq_B (E',R')$ and
$(E',R') \leq_B (E,R)$ hold. Two relations which are Borel bireducible
 to each other are said to  have the same complexity.  We write $(E,R)
 <_B (E',R')$ when
$(E,R) \leq_B (E',R')$ but $(E,R) \not \sim_B (E',R')$.

In the theory of classification of analytic equivalence relations on
 Polish spaces,
one tries to classify those relations up to Borel bireducibility.
Even for Borel relations, the situation is quite complicated, but
 there are a
number of natural milestones. They correspond
 to canonical equivalence relations on some classical Polish spaces. 

We thus have a scale of canonical relations, and given an equivalence
 relation
on a Polish space, we wish to locate it on this scale of complexity.

\vs 0.2cm

{\bf 1.2.  Complexity of analytic equivalence relations on Polish
 spaces.} We give some of these natural milestones and the relations
 between them.
The relation $(n,=)$ of equality on $n \in \N$ is the canonical
 example of a relation with $n$ classes. We also define $(\om,=)$ and
 $(2^\om,=)$. Because of their cardinalities, it is clear that
\vspace{-0,3cm}
 $$(1,=) <_B ... <_B (n,=) <_B (\om,=) <_B (2^\om,=). \vspace{-0,3cm}
 \leqno (1)$$
\hs 0.5cm The next relation is $(2^\om,E_0)$, or in short $E_0$. It is
 defined on $2^\om$ by
\vspace{-0,3cm}
$$\alpha E_0 \beta \Leftrightarrow \exists m \ \forall n \geq m, \
 \alpha(n)=\beta(n), \vspace{-0,3cm}$$
and it is well-known and not difficult to see that it satisfies 
\vspace{-0,3cm}
$$(2^\om,=) <_B E_{0}. \vspace{-0,3cm} \leqno (2)$$

By a theorem of Silver \cite{Si} and a theorem of
 Harrington-Kechris-Louveau
 \cite{HKL}, this list is extensive for those Borel equivalence
 relations
 which are Borel reducible to $E_0$. This is false for analytic
 equivalence
 relations, see \cite{R} for more details.

 After $E_0$, the order is no longer linear and the natural examples
 fall into one of two groups.

 A first family of milestones is given by Borel action of Polish
groups on Polish spaces. Given such an action of a Polish group $G$ on
 a
Polish space $X$, the orbit relation $E_G^X$ on $X$, defined by
 $x E_G^X y \Leftrightarrow \exists g \in G: y=g.x$, is an equivalence
 relation on $X$. It is called a {\em $G$-equivalence relation} and it
 is analytic.
It turns out that given a Polish group $G$, there is always a
$G$-equivalence relation which is $\leq_B$-maximum
among all possible $G$-equivalence relations on Polish spaces
 \cite{BK}. This
equivalence relation is denoted by $E_G^{\infty}$, without explicit
 reference
to the Polish space on which it is defined.

 Of particular interest are $E_{F_2}^{\infty}$, where 
$F_2$ is the free group with $2$ generators,
 $E_{S_{\infty}}^{\infty}$, where
$S_\infty$ is the group of permutations of the integers, and
 $E_{G_0}^{\infty}$, where $G_0$ is the group of homeomorphisms of the
 Hilbert cube.

In fact, $E_{F_2}^{\infty}$ is $\leq_B$-maximal among Borel
 equivalence
relations on Polish spaces for which each equivalence class is
 countable
or equivalently, among $G$-equivalence relations for countable groups
 $G$
\cite{FM}, \cite{JKL}.
 The relation  $E_{G_0}^{\infty}$ is $\leq_B$-maximum among all
 $G$-equivalence
 relations for Polish groups $G$ (Theorem 9.18 in \cite{Ke2} and
 Theorem
 2.3.5. in \cite{BK}). We have 
\vspace{-0,3cm}
$$E_0 <_B E_{F_2}^{\infty} <_B
E_{S_\infty}^{\infty} <_B  E_{G_0}^{\infty}. \vspace{-0,3cm} \leqno
 (3)$$

\

On the other hand, not all Borel equivalence relations are Borel
 reducible to equivalence relations associated to Borel actions of
 Polish groups. This is the second family of milestones "on the other
 side". 

The relation $E_1$ \cite{KL} is defined on ${\R}^{\om}$ by 
\vspace{-0.3cm}
$$\alpha E_1 \beta \Leftrightarrow \exists m \ \forall n \geq m, \
 \alpha(n)=\beta(n). \vspace{-0.3cm}$$
\hs 0.6cm It is not reducible to any $G$-equivalence relation for any
 Polish group $G$ \cite{KL}.

There exists a $\leq_B$-maximum equivalence
 relation $E_{K_\sigma}$ among $K_\sigma$ equivalence relations
 \cite{R2}. It is said to be {\em $K_\sigma$-complete}  and satisfies
\vspace{-0.3cm}
$$E_1 <_B E_{K_\sigma}.\vspace{-0.3cm} \leqno (4)$$
\hs 0.6cm  Rosendal \cite{R2} has found useful representations of this
 equivalence relation, we will use the following one, which was
 actually the starting point for this paper. Let $X_0$ be the set
$\Pi_{n \geq 1} n$.  The relation $H_0$ on $X_0$ is defined by
\vspace{-0,3cm}
$$\alpha H_0 \beta \Leftrightarrow \exists N \  \forall k,
|\alpha(k)-\beta(k)| \leq N. \vspace{-0,3cm}$$
\hs 0.6cm In  \cite{R2} it was proved that the relation $H_0$ is
$K_{\sigma}$-complete, that is, Borel bireducible with
 $E_{K_{\sigma}}$.

\

We may add canonical examples to our list using the operation $+$
 defined as
follows, see e.g. \cite{FS}. Let $E$ be an analytic equivalence
 relation on a Polish space $X$.
Then $E^+$ is the (also analytic) equivalence relation defined on
 $X^{\om}$ by
\vspace{-0,3cm}
$$(x_n) E^+ (y_n) \Leftrightarrow \forall n \ \exists m,p: (x_n E y_m)
 \wedge
(y_n E x_p).$$
\hs 0.6cm For example, $(2^\om,=)^+$, also written $=^{+}$,  is the
 relation
of equality of countable subsets of $2^{\om}$. By properties of
$E_{S_{\infty}}^{\infty}$ and the 'jump' properties of $+$, see
 \cite{FS}, or \cite{HIKO} (where $=^{+}$ is called $E_{countable}$), 
it satisfies
\vspace{-0,3cm}
$$E_{F_2}^{\infty}<_B =^{+} <_B E_{S_{\infty}}^{\infty}
 \vspace{-0,3cm} \leqno (5)$$
\hs 0.6cm Finally, it is also known that there exists a $\leq_B$-maximal element
 among analytic equivalence relations on Polish space, it is denoted
 by $E_{\Sigma_1^1}$ \cite{LR}.
Representations of this relation are for example isometric biembeddability
 between
separable Banach spaces or isometric biembeddability between metric Polish
 spaces
\cite{LR}.

The $\leq_B$-relations (1), (2), (3), (4)  and (5) can be summarized
 as follows:

\begin{figure}[htbp]
\begin{center}
\psfrag{a}{{\scriptsize {$(n, =)$}}}
\psfrag{b}{{\scriptsize {$(\omega, =)$}}}
\psfrag{c}{{\scriptsize {$(2^{\omega}, =)$}}}
\psfrag{d}{{\scriptsize {$E_0$}}}
\psfrag{e}{{\scriptsize {$E^{\infty}_{F_2}$}}}
\psfrag{f}{{\scriptsize {$=^{+}$}}}
\psfrag{h}{{\scriptsize {$E^{\infty}_{S_{\infty}}$}}}
\psfrag{i}{{\scriptsize {$E^{\infty}_{G_0}$}}}
\psfrag{j}{{\scriptsize {$E_{\Sigma^{1}_{1}}$}}}
\psfrag{k}{{\scriptsize {$E_{K_\sigma}$}}}
\psfrag{l}{{\scriptsize {$E_1$}}}
\psfrag{n}{{\scriptsize {$(1, =)$}}}
\psfrag{q}{{\scriptsize {\hbox{Figure 1: simplified diagram of
 complexity of analytic equivalence relations on Polish spaces.}}}}
\includegraphics[scale=0.30]{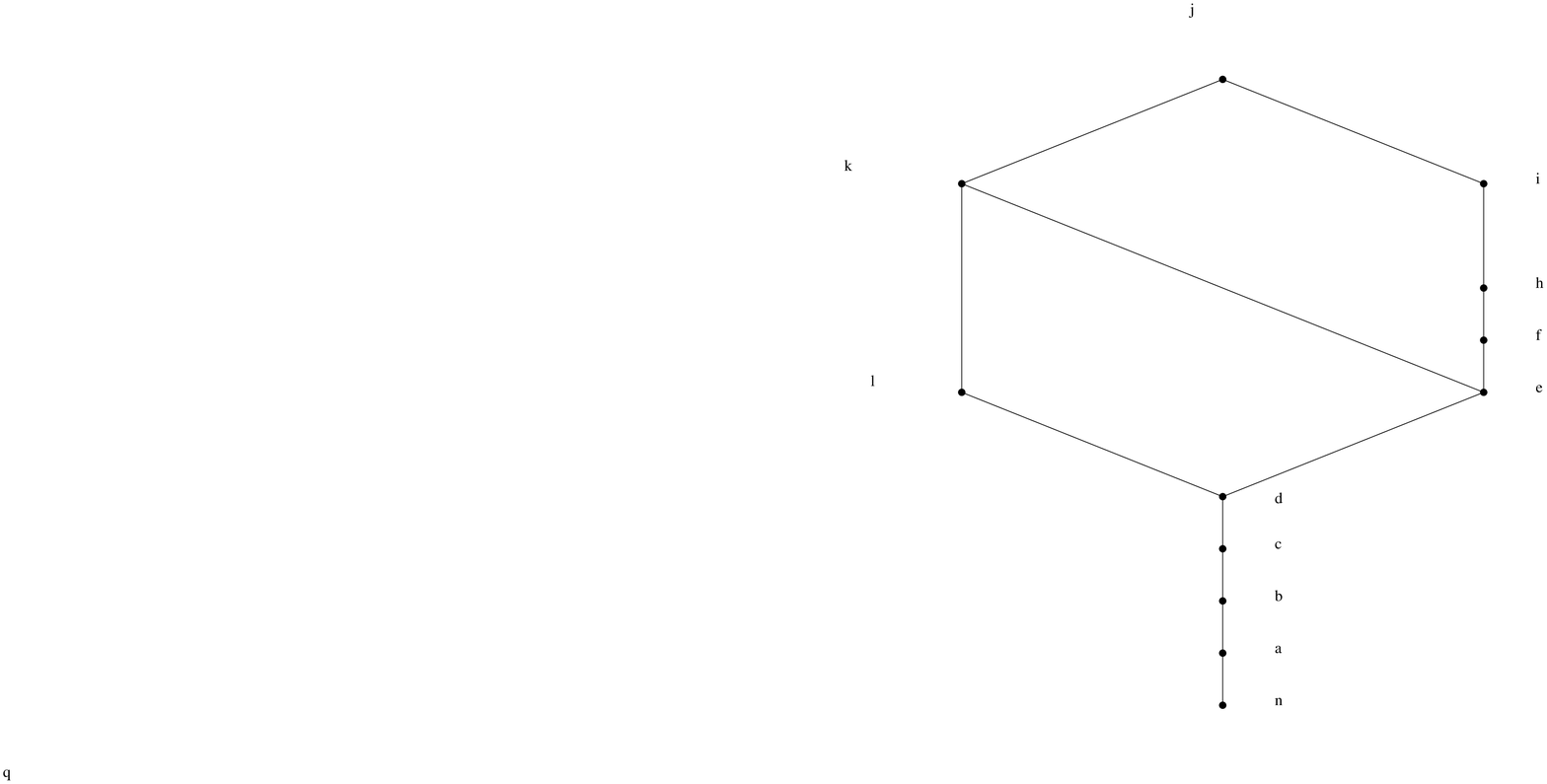}
\label{figura2004}
\end{center}
\end{figure}
\vspace{-0,5cm}
 {\bf 1.3. Complexity of isomorphism of separable Banach spaces.} Now
 if one is interested in Banach space theory, there are two possible
 directions in relation with the theory of classification of analytic
 equivalence relations up to Borel bireducibility.

The first one is to determinate the general $\leq_B$ complexity of
 isomorphism between separable Banach spaces, or at least a lower
 bound for this complexity. In other words, to show that isomorphism
 between separable Banach spaces reduces rather complex equivalence
 relations.

The question is also  asked for other natural equivalence
 relations of
interest in Banach space theory, such as Lipschitz isomorphism,
biembeddability, complemented biembeddability, isometry and so on.

 There is a natural setting for this question. All separable Banach
 spaces can be seen as subspaces of an isometrically universal Banach
 space such as $C([0,1])$. The set of such subspaces can be equipped
 with the Effros-Borel structure, see
e.g.  \cite{B},
which turns it into a standard Borel space.

 However in practice, 
one deals with particular examples, like subspaces generated by
 subsequences of a given basis, or block-subspaces of a given basis,
 and uses an ad-hoc topology 
associated
to the set of Banach spaces used for the reduction.
For example, all Schauder bases can be seen as subsequences of the
universal basis of Pe\l czy\'nski \cite{LT}, and thus a Banach space with
 a basis can be represented as an equivalence class on $2^\om$.

There are only a few results in that direction, and they are recent. 
 B. Bossard \cite{B} proved that isomorphism between Banach spaces is
 analytic
  non-Borel, and Borel reduces $E_0$. Using Tsirelson's space,
 Rosendal \cite{R} improved the result to $E_1$, which implies that
 isomorphism between separable Banach spaces is not associated to
  the Borel action of a Polish group on a Polish space.

The other direction of research is to try to relate the complexity of
 the
 relation of isomorphism between subspaces of a given separable Banach
 space
 $X$ to geometrical properties of $X$.

 Indeed, by the solution to the
 homogeneous Banach space problem, given by W. T. Gowers \cite{Go1}
 and R. A. Komorowski, N. Tomczak-Jaegermann \cite{TJ}, if a Banach
 space $X$ has only one class of isomorphism of subspaces, then $X$
 must be isomorphic to $l_2$.
But it is not known if, for example, there exists a Banach space, other
 than
 $l_2$, with at most $2$,
 or even at most $\omega$ classes of isomorphism of subspaces.
 The following question, asked by G. Godefroy, is then natural: if the
 complexity of isomorphism between subspace of $X$ is low (in the
 sense of
 cardinality or more interestingly in the sense of $\leq_B$), then
 what
 geometrical, regularity properties must $X$ satisfy?

 It turns out that a first natural threshold for this question is the
 relation $E_0$.
Indeed, recently the first named author and Rosendal \cite{FR2} have
 defined
 a Banach space $X$ to be {\em ergodic} if $E_0$ is Borel reducible to
 isomorphism between subspaces of $X$ and
have obtained various results about ergodic spaces.

 Rosendal \cite{R1} noticed that
 hereditarily indecomposable Banach spaces are ergodic, and also proved that
 an unconditional basis of a non-ergodic Banach space must have a subsequence
 such that all further subsequences span isomorphic susbpaces. In
 \cite{FR1} and \cite{R1} it is proved that a non-ergodic Banach space
 $X$ with an
 unconditional basis must be isomorphic to its square, its hyperplanes
 and
 more generally to $X \oplus Y$ for any subspace  $Y$ generated by a
 subsequence of
 the basis. This was already obtained by Kalton \cite{K2} for spaces
 with an
 unconditional basis and only countably many classes.
Finally,
 in \cite{FR2} it is
 proved that a non-ergodic Banach space must contain a subspace $X$
 with an unconditional basis which is isomorphic to $X \oplus Y$ for
 all block-subspaces $Y$ of $X$.

\

In this paper, we work with the classical Banach spaces $c_0$, $l_p$
 and $L_p$, $1 \leq p <2$, to  provide new results concerning the two
 directions of research.

\vs 0.2cm

{\bf 1.4. Organization of the paper.} In the next section we Borel
 reduce the relation $E_{K_{\sigma}}$ to
isomorphism between 
subspaces of $l_p, 1 \leq p<2$ (Theorem 2.6).  Our main tool for this
 is a theorem of P. G. Casazza and  N. Kalton about uniqueness of
 unconditional structure for
Banach spaces, see Theorem 2.3. This result will allow us to prove
 that certain spaces with
unconditional bases are isomorphic if and only if their canonical
 bases are equivalent. It was Kalton who suggested that this paper
 contained an answer to the problem of the number of non-isomorphic
 subspaces of $l_p$.

In the third section we also reduce the relation $E_{K_{\sigma}}$ to 
isomorphism between 
subspaces of $c_0$ (Theorem 3.3). There we use another theorem of  P.
 G. Casazza and
N. Kalton (Theorem 3.1) which is a strengthening of Theorem 2.3 in the case of
$c_0$-sums. In
 particular,
 our results show that $c_0$ and $l_p,
 1 \leq p<2$ are ergodic. This extends to Banach spaces with an unconditional basis
 with the shift property and satisfying a lower $p$-estimate, $1 \leq p <2$ (Theorem 2.7). 

In the fourth section, using more simple techniques, we reduce the
 relation $=^{+}$ to
 isomorphism between subspaces of $L_p$, $1 \leq p <2$ (Theorem 4.1).
In combination with Theorem 2.6, we deduce that $E_{K_{\sigma}} \otimes =^+$  
 is Borel reducible to isomorphism between subspaces of $L_p, 1 \leq p
 <2$. 
 Thus isomorphism between separable Banach spaces is
 not
 reducible to the equivalence relation associated to the Borel action
 of a
 Polish group, but reduces $G$-equivalence relations for
 non-trivial actions of such groups $G$. 

Finally, in the last section  we  present a diagram (figure 2)
 containing known facts about complexity  of isomorphism between
 subspaces  of a Banach space. Then,  we  point out  some open problems and a
 conjecture  related to our results.

\vs 0.2cm

{\bf 1.5. Notation.} We shall write $X \simeq Y$ to mean that two
 Banach spaces $X$ and $Y$ are isomorphic, $X \st{L} {\simeq} Y$ to mean that
 they are Lipschitz isomorphic, $X \st{c}{\hr} Y$ to mean that $X$ is
 isomorphic to a complemented subspace of $Y$, and $X \st{c} {\simeq} Y$ to
 mean that $X \st{c}{\hr} Y$ and $Y \st{c}{\hr} X$. For Banach spaces
 $X$ and $Y$ written as $c_0$-sums or $l_p$-sums of $l_q^n$ spaces, we
 shall abuse notation by writing $X \sim Y$ to mean that the canonical
 bases are equivalent.

If $(X_n)_{n \in \N}$ is a sequence of Banach spaces, we will denote
 the $l_p$-sum  $(\sum_{n \in \N} \oplus X_n)_{l_p}$
of the $X_n$'s by $l_p(X_n)_{n \in \N}$. The $l_p$-sum of infinitely
 many copies of a Banach space $X$ is denoted as usual $l_p(X)$.
If the Banach spaces $X_{n}$, $n \in \N$, are given with canonical
 bases, then $l_p(X_n)_{n \in \N}$ has a corresponding canonical basis
associated to a given bijection between $\N$ and $\N^2$.
We use the similar notation for $c_0$-sums.

\vs 0.4cm

{\bf 2. Reductions of $E_{K_\sigma}$ to isomorphism between subspaces
 of $l_p$, $1 \leq p <2$.}

\vs 0.2cm

 Rosendal \cite{R1} proved that $E_{K_{\sigma}}$ is Borel reducible to
 equivalence
between Schauder basis in Banach spaces. We start by rephrasing his
 proof in a
slightly higher generality (Lemma 2.1 and Corollary 2.2).  

\vs 0.3cm

{\bf Lemma 2.1:}  {\it Let $(K_n)_{n \in \N}$ be a sequence of
 integers,
  $(p_n)_{n \in \N}$, $(q_n)_{n \in \N}$  be bounded sequences of
 reals larger than $1$ and  $1 \leq p < +\infty$. Then 
\vspace{-0,2cm}  
$$l_p(l_{p_n}^{K_n})_{n \in \N} \sim
l_p(l_{q_n}^{K_n})_{n \in \N} \Leftrightarrow \exists C>0, \forall n
 \in \N, |p_n-q_n| \leq \frac{C}{\log K_n}, \vspace{-0,3cm}$$
and the similar result is valid for $c_0$-sums.}

\vs 0.2cm

{\it Proof.} By a classical consequence of H\"older's inequality, the
constant of equivalence $c_n$ between the canonical bases of
$l^{K_n}_{p_n}$ and
$l^{K_n}_{q_n}$ is
$K_n^{|1/p_n-1/q_n|}$ \cite{T2}.

Let $c$ be an upper bound for the sequences $(p_n)_{n \in \N}$ and
 $(q_n)_{n \in \N}$, then
\vspace{-0,1cm}
$$e^{|p_n-q_n|\log K_n/c^2} \leq c_n \leq
e^{|p_n-q_n| \log K_n}. \vspace{-0,1cm}$$

It follows that if the canonical bases of $l_p(l^{K_n}_{p_n})_{n \in
 \N}$ and $l_p(l^{K_n}_{q_n})_{n \in \N}$ (resp. $c_0(l^{K_n}_{p_n})_{n \in
 \N}$ and $c_0(l^{K_n}_{q_n})_{n \in \N}$)
are $C$-equi\-valent, then
for all $n$, $|p_n-q_n| \leq c^2\log C/{\log K_n}$. 

Conver\-sely, if
for all $n$, $|p_n-q_n| \leq M/{\log K_n}$, then the ca\-nonical bases
 of $l_p(l^{K_n}_{p_n})_{n \in \N}$ and 
$l_p(l^{K_n}_{q_n})_{n \in \N}$ (resp. $c_0(l^{K_n}_{p_n})_{n \in \N}$ and 
$c_0(l^{K_n}_{q_n})_{n \in \N}$) are $(e^M)^2$-equivalent. $\cqd$ 

\vs 0.3cm

Let $(K_n)_{n \in \N}$  be a sequence of integers, $(p_n)_{n \in \N}$
 be a sequence of real numbers greater than $1$ and $1 \leq p
 <+\infty$. We recall that $X_0$ denotes the set
$\Pi_{n \geq 1} n$, and that the  relation $H_0$ on $X_0$ is defined by
\vspace{-0,3cm}
$$\alpha H_0 \beta \Leftrightarrow \exists N \  \forall k,
|\alpha(k)-\beta(k)| \leq N. \vspace{-0,3cm}$$

  For $\al \in X_0$, we denote by
$l_p(l_{p_n}^{K_n}(\alpha))$ the Banach space
$$l_p(l_{p_n}^{K_n}(\alpha))=
(\sum_{n} \oplus
l^{K_n}_{p_n+\frac{\alpha(n)}{\log K_n}})_{l_p}, \vspace{-0,3cm}$$
and we use the similar definition for $c_0$-sums.

\vs 0.3cm

{\bf Corollary 2.2:} {\it Let $(K_n)_{n \in \N}$ be a sequence of
 integers,
  $(p_n)_{n \in \N}$  be a sequence of reals larger than $1$ such that
$(p_n+\frac{n}{\log K_n})_{n \in \N}$ is bounded and 
 $1 \leq p < +\infty$. Then for all $\alpha$ and $\beta$ in $X_0$,
\vspace{-0,3cm} 
$$\alpha H_0 \beta \Leftrightarrow  l_p(l_{p_n}^{K_n}(\alpha)) \sim
l_p(l_{p_n}^{K_n}(\beta)),$$ and
  the similar result is valid for $c_0$-sums.}

\hs 0.6cm

 It is known that for $1 \leq p \leq r \leq 2$ and $\e>0$, $l_p$ contains
 $1+\e$-isomorphic copies of $l_r^n$, in fact $L_r$
 is isometric to some subspace of $L_p$, see e.g. \cite{LT}. It
 follows that for any sequence $(p_n)_{n \in \N}$ of reals such that
for all $n$, $p \leq p_n \leq 2$ and any sequence of integers
 $(K_n)_{n \in \N}$, the space
$l_p(l_{p_n}^{K_n})_{n \in \N}$ is isomorphic to a subspace of $l_p$.

Our main ingredient will be the following theorem of Casazza and
 Kalton
\cite{CK1}, which can be thought of as a first step towards uniqueness
 of unconditional
structure for Banach spaces which are sufficiently far from $l_2$. 

\

We refer to \cite{LT}, \cite{K} for the definition of and background about Banach lattices.
If $X$ and $Y$ are Banach lattices, a bounded linear operator $V:X \rightarrow
Y$ is called a {\em lattice homomorphism} if $V(x_1 \vee x_2)=Vx_1 \vee Vx_2$
for all $x_1,x_2 \in X$.
 Following \cite{CK1}, define a Banach lattice $X$ to
be {\em sufficiently lattice-euclidean} if there exists $C \geq 1$ such thar
for all $n \in \N$, there exist operators $S:X \rightarrow l_2^n$ and $T:l_2^n
\rightarrow X$ such that $ST=I_{l_2^n}$, $\norm{S}\norm{T} \leq C$ and such
that
$S$ is a lattice homomorphism.
This is equivalent to saying that $l_2$ is finitely representable as a 
complemented sublattice of $X$.

A Banach space with a $1$-unconditional basis $(x_n)_{n \in \N}$ is naturally
considered as a Banach lattice by defining $$\sum_{n \in \N}a_n x_n \geq 0
\Leftrightarrow \forall n \in \N, a_n \geq 0.$$
It is classical to consider a Banach space $X$ with a $C$-unconditional basis,
$C \geq 1$, as a
Banach lattice as well, with the same definition of $\leq$ and with the
restriction of having to add a constant in some inequalities; alternatively,
one may equip $X$ with an equivalent norm which turns $(x_n)_{n \in \N}$
$1$-unconditional, and the results concerning the Banach lattice structure of
$X$ can be transfered back to the initial norm,
up to to the constant of equivalence.

For an unconditional basis $(x_n)_{n \in \N}$ of a Banach space (seen as a
Banach lattice), being
sufficiently lattice-euclidean is the same as
 having, for some $C \geq 1$ and every $n \in \N$,
a $C$-complemented, $C$-isomorphic copy of $l_2^n$
 whose basis is disjointly supported on $(x_n)_{n \in \N}$. 

\vs 0.3cm

{\bf Theorem 2.3: (Casazza -Kalton \cite{CK1}}) {\it Let $X$ be  a Banach
  space with an unconditional basis and $(y_n)_{n \in \N}$ be an
 unconditional basic sequence in $X$ which is not sufficiently
 lattice-euclidean and spans a complemented subspace of $X$. Then
 $(y_n)_{n
 \in \N}$ is equivalent to a
 sequence of disjointly supported vectors which spans a complemented
 subspace in $X^N$ for some $N$.} 

\vs 0.3cm

 Let $X$ be a Banach space with a Schauder decomposition
 $X=\sum_{i=1}^{+\infty}
 \oplus X_i$. We shall say that vectors $x$ and $y$ in $X$ are {\em
 successive} and write $x<y$ if there exists intervals of integers $E$
 and $F$
 such that $\max(E)<\min(F)$, $x \in \sum_{i
 \in E} X_i$, and $y \in \sum_{i \in F} X_i$. 

We say that
 the Schauder decomposition of $X$ 
satisfies a {\em lower $p$-estimate with constant $C \geq 1$} if for
 any successive
 vectors $x_1<\cdots<x_n$ in $X$,
 $(\sum_{i=1}^n \norm{x_i}^p)^{1/p} \leq C \norm{\sum_{i=1}^n x_i}$.
For $1 \leq p \leq +\infty$, the conjugate $p'$ of $p$ is as usual
 defined by 
$\frac{1}{p}+\frac{1}{p'}=1$ (with $\frac{1}{+\infty}=0$).

\vs 0.3cm

{\bf Lemma 2.4:} {\it Fix $1 \leq p <2$, $(K_n)_{n \in \N}$ a sequence
 of integers and   $(p_n)_{n \in \N}$ 
a sequence of real numbers which is bounded below by $p$. 
Let $r=\sup_{n \in \N} p_n$ and  $r^{\prime}$ be the conjugate of $r$.
 Suppose
$X=\sum_{n=1}^{+\infty} \oplus l_{p_n}^{K_n}$ is a Schauder
 decomposition of $X$ satisfying a lower $p$-estimate with constant $C
 \geq 1$. Then for all $k \in \N$, for all vectors
$y_1,\ldots,y_k$ in  $X$ which are disjointly supported on its
 canonical basis,
\vspace{-0,3cm}
$$\sum_{i=1}^k \norm{y_i} \leq C k^{1/r^{\prime}} \norm{\sum_{i=1}^k
 y_i}. \vspace{-0,3cm}$$}

\vs 0.2cm

{\it Proof.} We may assume that $r<+\infty$. Let $y_1,\ldots,y_k$ be
 as above. For each $1 \leq i \leq k$, we write
$y_i=\sum_{n=1}^{+\infty} y_{in}$, where $y_{in}$ is the projection of
 $y_i$ onto the $l_{p_n}^{K_n}$ summand. Then
\vspace{-0,3cm}
$$C\norm{\sum_{i=1}^k y_i}=C\norm{\sum_{n=1}^{+\infty} \sum_{i=1}^k
 y_{in}}
\geq (\sum_{n=1}^{+\infty} (\sum_{i=1}^{k}
 \norm{y_{in}}^{p_n})^{\frac{p}{p_n}})^{\frac{1}{p}}.
 \vspace{-0,3cm}$$

Denote $a_{in}=\norm{y_{in}}^{p}$ and $\alpha=r/p$, $\alpha_n=p_n/p$.
Then
\vspace{-0,3cm}
$$C^{p}\norm{\sum_{i=1}^k y_i}^p \geq  \sum_{n=1}^{+\infty} (\sum_{i=1}^{k} a_{in}^{\alpha_n})^{\frac{1}{\alpha_n}} 
\geq  \sum_{n=1}^{+\infty} k^{-1/{\alpha_n}^{\prime}} \sum_{i=1}^k a_{in}, \vspace{-0,3cm}$$
by H\"older's inequality, since $\alpha_n \geq 1$.
Now for every $n \in \N$, $\alpha \geq \alpha_n$, so
\vspace{-0,3cm}
$$C^{p}\norm{\sum_{i=1}^k y_i}^p \geq  k^{-1/{\alpha}^{\prime}} \sum_{n=1}^{+\infty} \sum_{i=1}^k a_{in}. \vspace{-0,3cm}$$
\hs 0.6cm On the other hand,
\vspace{-0,3cm}
$$(\sum_{i=1}^{k} \norm{y_i})^p=(\sum_{i=1}^k (\sum_{n=1}^{+\infty} a_{in})^{\frac{1}{p}})^p \leq k^{\frac{p}{p^{\prime}}}
\sum_{i=1}^{k} (\sum_{n=1}^{+\infty} a_{in}), \vspace{-0,3cm}$$
once again by H\"older's inequality, since $p \geq 1$.
Finally,
\vspace{-0,3cm}
$$(\sum_{i=1}^k \norm{y_i})^p \leq C^p k^{\frac{p}{p^{\prime}}+
\frac{1}{\alpha^{\prime}}} \norm{\sum_{i=1}^{k}y_i}^p,  \vspace{-0,3cm}$$
so
\vspace{-0,3cm}
$$\sum_{i=1}^k \norm{y_i} \leq C k^{\frac{1}{p^{\prime}}+
\frac{1}{p\alpha^{\prime}}} \norm{\sum_{i=1}^{k}y_i}, \vspace{-0,3cm}$$
and the fact that 
$\frac{1}{p^{\prime}}+
\frac{1}{p\alpha^{\prime}}=\frac{1}{r^{\prime}}$ concludes the proof. $\cqd$ 

\vs 0.3cm

To  prove the next proposition we need to recall that two unconditional sequences $(u_n)_{n \in \N}$ and $(v_n)_{n \in \N}$ in a Banach space $X$ are said to be {\it permutatively equivalent} if there is a permutation $\pi$ of $\N$ so that $(u_n)_{n \in \N}$ and $(v_{\pi (n)})_{n \in \N}$ are equivalent.

\vs 0.3cm

{\bf Proposition 2.5:} {\it Let $(K_n)_{n \in \N}$ be a sequence of integers, $(p_n)_{n \in \N}$, $(q_n)_{n \in \N}$  sequences of real numbers and $1 \leq p <2$. Assume}

(1) {\it $p<p_n<2$ and $p<q_n<2$, for all $n \in \N$;}

(2) {\it $(p_n)_{n \in \N}$ and  $(q_n)_{n \in \N}$ are decreasing sequences:}

(3) {\it $K_1 \geq 4$ and  $K_n \geq n^2 K_{n-1}$, for all $n \geq 2$.

Then whenever $l_p(l_{p_n}^{K_n})_{n \in \N} \st{c}{\hr} l_p(l_{q_n}^{K_n})_{n
  \in \N} \oplus F$, for some finite-dimensional space $F$,  
there exists $C>0$ such that  $p_n-q_n \leq C/\log K_n$, for all $n \in \N$.}

\vs 0.2cm

{\it Proof.} Note that by Lemma 2.4, any
 dis\-joint\-ly sup\-por\-ted se\-quen\-ce of vectors $x_1,\ldots,x_k$ 
 in $l_p(l_{p_n}^{K_n})_{n \in \N}$ satis\-fies 
\vspace{-0,3cm}
$$\sum_{i=1}^{k} \norm{x_i} \leq k^{1/q_1^{\prime}} \norm{\sum_{i=1}^k x_i}, \vspace{-0,3cm}$$
and ${q_1}^{\prime} >2$. So $l_p(l_{p_n}^{K_n})_{n \in \N}$ is not sufficiently lattice-euclidean. 
By Theorem 2.3, for some $N$, the canonical basis of
 $l_p(l_{p_n}^{K_n})_{n \in \N}$ is $C$-equivalent to a disjointly supported
 sequence in $(l_p(l_{q_n}^{K_n})_{n \in \N} \oplus F)^N$. Modifying $N$ and
 $C$ we
 may assume $F=\{0\}$. Then
without loss of generality we may write this space as $l_p(l_{q_n}^{NK_n})_{n \in \N}$ (the canonical bases are permutatively equivalent).
Take $k \geq k(N)$, where $k(N)$ is such that this condition ensures 
$\frac{K_k}{2} \geq \sum_{i=1}^{k-1} NK_i$; it exists by condition (3).
The canonical basis of $l_{p_k}^{K_k}$ is $C$-equivalent to a
disjointly supported sequence in $l_p(l_{q_n}^{NK_n})_{n \in \N}$. By the condition on $K_k$,
we see that  the canonical basis $(e_i)_{i \in \N}$ of $l_{p_k}^{K_k/2}$ is $C$-equivalent to a
disjointly supported sequence $(f_i)_{i \in \N}$ in $(\sum_{n \geq k} \oplus l_{q_n}^{NK_n})_{l_p}$. 
We may now apply Lemma 2.4 to the sequence $(f_i)_{i \in \N}$.
 As $q_k=\max \{q_n, n\geq k \}$,
\vspace{-0,3cm}
$$C(K_k/2)^{1/p_k} \geq \norm{\sum_{i=1}^{K_k/2} f_i} \geq 
(K_k/2)^{-1/{q_k}^{\prime}} 
\sum_{i=1}^{K_k/2} \norm{f_i} \geq (K_k/2)^{1/q_k}/C. \vspace{-0,3cm}$$
\hs 0.6cm Consequently
$$(K_k/2)^{1/q_k-1/p_k} \leq C^2, \vspace{-0,3cm}$$
and
\vspace{-0,3cm}
$$p_k-q_k \leq 4(1/q_k-1/p_k) \leq 8\log C/\log(K_k/2) \leq 16\log C/\log K_k. \vspace{-0,3cm}$$
\hs 0.6cm This is true for any $k \geq k(N)$, so the proposition is proved. $\cqd$ 

\vs 0.3cm

{\bf Theorem 2.6:} {\it Suppose $1 \leq p <2$. Then the relation $E_{K_{\sigma}}$ is Borel reducible to isomorphism, to Lipschitz isomorphism and to complemented biembeddability between
subspaces of $l_p$. Indeed, there exist a sequence of integers
$(K_n)_{n \in \N}$, a sequence of reals $(p_n)_{n \in \N}$ with $p<p_n<2$ for all $n$, such
that the following are equivalent for
 all $\alpha$ and $\beta$ in $X_0$:}

(1)  {\it $\alpha H_0 \beta$.} 

(2) {\it  $l_p(l_{p_n}^{K_n}(\alpha)) \sim l_p(l_{p_n}^{K_n}(\beta))$.}

(3)  {\it $l_p(l_{p_n}^{K_n}(\alpha)) \simeq l_p(l_{p_n}^{K_n}(\beta))$.}

(4)  {\it $l_p(l_{p_n}^{K_n}(\alpha)) \st{L} {\simeq}  l_p(l_{p_n}^{K_n}(\beta))$.}

(5)  {\it $l_p(l_{p_n}^{K_n}(\alpha)) \st{c} {\simeq} l_p(l_{p_n}^{K_n}(\beta))$.}

\vs 0.2cm

{\it Proof.} We choose $(K_n)_{n \in \N}$ satisfying (3) of Proposition 2.5, and $(p_n)_{n \in \N}$ such that $p_1+1/\log K_1<2$,  $p<p_n<2$ for all $n$, and $\frac{n+1}{\log K_{n+1}} \leq p_n-p_{n+1}$. This is certainly possible if $\sum_{n+1}^{+\infty} \frac{n}{\log K_n}$ is small enough.
Then it is clear that the conditions of Proposition 2.5 are satisfied for any 
two sequences $(p_n+\frac{\alpha(n)}{\log K_n})_{n \in \N}$
and $(p_n+\frac{\beta(n)}{\log K_n})_{n \in \N}$. It follows that (5) implies (1).
That (4) implies (5), that is,  Lispchitz isomorphism implies complemented
biembeddability, comes from the fact that the spaces
considered are separable dual spaces (Theorem 2.4 in \cite{HM}). 
(1) implies (2) by Lemma 2.1 and Corollary 2.2 and the rest is obvious.
\pff $\cqd$ 

\vs 0.3cm

Using a similar proof as in the previous theorem we get the following result. An unconditional basis for
a Banach space $X$ is said to have the shift property if any normalized block-sequence 
$(x_n)_{n \in \N}$
in $X$ is equivalent to $(x_{n+1})_{n \in \N}$.
\vs 0.2cm

{\bf Theorem 2.7:} {\it Let $X$ be a Banach space  with an unconditional basis with the shift property,
which satisfies a lower
  $p$-estimate for some $1 \leq p<2$. Then $X$ is ergodic.}

\vs 0.2cm

{\it Proof.} Let $X$ be such a space.  By Krivine's theorem (see e.g. \cite{M}), $l_r$ is block-finitely represented in $X$
for some $r$, and by the lower estimate,  we have that $r \leq p$. But then all $l_{r'}$ for $r \leq r'
\leq 2$, and in particular $p \leq r'\leq 2$ are
finitely represented in $X$ (with constant $2$ say). We may then associate to each $\alpha \in X_0$ a
subspace $X(\alpha)$ of $X$ which is a direct sum on the basis of
$l_{p_n}^{K_n}$s for some $p_n$s in $]p,2[$  as previously. The canonical
Schauder decomposition of the space
$X(\alpha)$ satisfies a lower $p$-estimate and is unconditional.
Also each $l_{p_n}^{K_n}$ has a canonical $1$-unconditional basis and $X$ satisfies the shift property,
so by \cite{CK0}, Proposition 2.3, the canonical basis of $X(\alpha)$ is unconditional.
 We may thus follow
the proof of Proposition 2.5 (note that if
$X=\sum
\oplus l_{q_n}^{K_n}$ satisfies a lower $p$-estimate with constant $C$, then $X^N \simeq
\sum \oplus l_{q_n}^{NK_n}$ satisfies a lower $p$-estimate as well, with a
constant depending on $C$, $N$ and $p$). We get that
\vspace{-0,2cm}
 $$X(\alpha) \simeq X(\beta) \Rightarrow \alpha
H_0 \beta. \vspace{-0,2cm}
$$ 
\hs 0.6cm It doesn't seem easy to get the converse without more
 information on the norm on $X$. We shall in fact reduce $E_0$ instead of
 $E_{K_{\sigma}}$. For this, consider $2^{\om}$ as a subset of $X_0$ by 
$j((\alpha(n))_{n \in \N})=(0,\alpha(1),2\alpha(2),3\alpha(3),\ldots)$. Then
 clearly, for any $\alpha, \ \beta \in 2^{\om}$, $\alpha E_0 \beta$ if and only if $j(\alpha) H_0 j(\beta)$, so from the above,
\vspace{-0,3cm} 
$$X(j(\alpha)) \simeq X(j(\beta)) \Rightarrow \alpha E_0 \beta. \vspace{-0,3cm}
$$ 
\hs 0.6cm But we also have 
\vspace{-0,3cm}
$$\alpha E_0 \beta \Rightarrow X(j(\alpha)) \simeq X(j(\beta)),\vspace{-0,3cm}$$ because if $\alpha E_0 \beta$ then $X(j(\alpha))$ and
$X(j(\beta))$ have canonical bases which differ by only a finite number of vectors. So $E_0$ is Borel reducible to isomorphism between subspaces of $X$. $\cqd$

\vs 0.4cm

{\bf 3. Reductions of $E_{K_\sigma}$ to isomorphism between subspaces of $c_0$.}

\vs 0.2cm

We now turn our attention to spaces of the form $c_0(l^{K_n}_{p_n}(\alpha))$. 
Note that as a $c_0$-sum of finite-dimensional spaces, every such space is isomorphic to a subspace of $c_0$. The previous results
concerning isomorphism and complemented biembeddability between 
subspaces of $l_p$, $1 \leq p <2$ extend by duality to quotients of $l_p, p >2$ and $c_0$, and thus by a classical
 theorem (Theorem 2.f.6 in \cite{LT}), also to subspaces of $c_0$. However, we
 shall improve these results by also reducing $E_{K_\sigma}$ to
complemented embeddability between subspaces of $c_0$.

 We recall that the definition of $\leq_B$ still makes
 sense when the relation is not an equivalence relation. In particular,  the
 $\leq_B$-classification of quasi-orders has consequences in the
 $\leq_B$-classification of equivalence relations, see \cite{LR}.

\vs 0.3cm

{\bf Theorem 3.1: (Casazza-Kalton \cite{CK2})} {\it Let $(K_n)_{n \in \N}$ be a sequence of integers and $(q_n)_{n \in \N}$  a decreasing sequence  of
 reals converging to $1$. Then  any unconditional basis of a complemented subspace of $c_0(l^{K_n}_{q_n})_{n \in \N}$
is permutatively equivalent to the canonical basis of
$c_0(l^{M_n}_{q_n})_{n \in \N}$, for some sequence $(M_n)_{n \in \N}$
such that $M_n/K_n$ is bounded.}

\vs 0.2cm

Since it is only implicit in their paper, we sketch how this theorem follows
from their results. We also refer to their paper for some definitions which
we would not use afterwards.
 
\vs 0.2cm

{\it Proof of Theorem 3.1.} Let $(u_k)_{k \in \N}$ be
an unconditional basis of a complemented subspace 
of $c_0(l^{K_n}_{q_n})_{n \in \N}$.
 According to  \cite{CK2} Corollary 2.5 and \cite{CK2} Theorem 1.1 we may assume that the $(u_k)$'s
are disjointly supported. By \cite{CK2} Theorem 3.2 we may assume that for all $n$ and $k$,
$\norm{u_k(n)}_{l^{K_n}_{q_n}}=0$ or $1$, and that there exists a
partition $\N=\cup_{n \in A} B_n$ of $\N$, such that the space spanned by $(u_k)_{k \in \N}$ is a $c_0$-sum of the spaces spanned by $(u_k)_{k \in B_n}$, and a $C$ such that for each $n$,
$(u_k)_{k \in B_n}$ is $C$-complemented, $C$-tempered (see the definition in \cite{CK2}). 
By the Claim in  \cite{CK2} Theorem 3.4, we see that for some $K$,
each $(u_k)_{k \in B_n}$ must be $K$ permutatively equivalent to
 the canonical basis of
$(\sum_{k \in D_n} \oplus
l^{P_k}_{q_k})_{c_0}$, with $|D_n| \leq K$ and
$P_k/K_k \leq K$. Furthermore, by (3) of  \cite{CK2} Theorem 3.2, for any $k$, the number of $n$'s such that $k \in D_n$ is uniformly bounded. The theorem follows.
\pff $\cqd$ 

\vs 0.3cm

{\bf Proposition 3.2:} {\it Let $(K_n)_{n \in \N}$ be a sequence of integers
  and $(p_n)_{n \in \N}$, $(q_n)_{n \in \N}$ be sequences of reals. Assume}

(1) {\it The sequences $(p_n)_{n \in \N}$ and $(q_n)_{n \in \N}$ are decreasing to $1$;}

(2) {\it $1<p_n+\frac{n}{\log K_n} <2$ and $1<q_n+\frac{n}{\log K_n} <2$, for all $n \in \N$;}

(3) {\it The sequence $(\frac{n}{\log K_n})_{n \in \N}$ is decreasing;}

(4) {\it $|q_n-p_m| \geq \min(m,n)/\log K_{\min(m,n)}$,  for all $m \neq n.$

Then whenever $c_0(l_{p_n}^{K_n})_{n \in \N}  \st{c}{\hr}           c_0(l_{q_n}^{K_n})_{n \in \N}$, it follows that 
$c_0(l_{p_n}^{K_n})_{n \in \N} \sim c_0(l_{q_n}^{K_n})_{n \in \N}.$}

\vs 0.2cm

{\it Proof.} By  Theorem 3.1, the
canonical basis of $c_0(l_{p_n}^{K_n})_{n \in \N}$, being equi\-valent to an uncon\-ditional ba\-sis
of a comple\-mented sub\-space of $c_0(l_{q_n}^{K_n})_{n \in \N}$, must be $C$-permu\-tatively equi\-valent to the
cano\-nical basis of a space $c_0(l^{M_n}_{q_n})_{n \in \N}$, for some constant $C$ and some sequence $(M_n)_{n \in \N}$
of integers.
 
Let us now fix $n \in \N$. Thus by the above, and using the symmetry
of the canonical bases of spaces $l_{p}$, 
there exists for $i=1,\ldots,k$ integers $A_i$, such that
$K_n=\sum_{i=1}^k A_i$, and an increasing sequence of
integers $n_i, i
\leq k$, such that the canonical basis of $l^{K_n}_{p_n}$ is
$C$-equivalent to the canonical basis of $(\sum_{1 \leq i \leq k} \oplus
l^{A_i}_{q_{n_i}} )_{c_0}$.

We first note that in particular the canonical basis of $l_{\infty}^k$ is $C$-equivalent to
the canonical basis of $l^k_{p_n}$, from which it follows that
\vspace{-0,3cm}
$$k^{1/p_n} \leq C, \vspace{-0,3cm}$$
therefore $k \leq C^{p_n} \leq C^2$.

It follows that $\max_{i} A_i \geq K_n/C^2$. Let $i$ be an integer where this 
maximum is attained and let $N=n_i$. From the fact that the canonical bases of  
$l^{A_i}_{p_n}$
and $l^{A_i}_{q_N}$ are $C$-equivalent, it follows that 
\vspace{-0,3cm}
$$A_i
^{|\frac{1}{p_n}-\frac{1}{q_N}|}
\leq C, \vspace{-0,3cm}$$
so by condition  (2),
 $$1/4 |p_n-q_N| \log A_i
\leq \log C.$$

We now assume that
 $n > 8\log C$, it follows from (2) that $\log K_n > 4\log C$,
and therefore 
\vspace{-0,3cm}
$$|p_n-q_N|  \leq
4\log C/(\log K_n - 2\log C) \leq 8\log C/\log K_n. \vspace{-0,3cm}$$ 

Now if $N \neq n$, then by conditions (3) and (4),
\vspace{-0,3cm} 
 $$|p_n-q_N| \geq \min(n,N)/\log
K_{\min(n,N)} \geq n/\log K_n. \vspace{-0,3cm}$$ 
\hs 0.6cm But this contradicts the
assumption that $n > 8\log C$.
It follows that $N=n$. In particular in the previous
inequality, we get
\vspace{-0,3cm}
$$|p_n-q_n|  \leq
 8\log C/\log K_n. \vspace{-0,3cm}$$

Finally, for all $n>8 \log C$, $|p_n-q_n| \leq 8\log C/\log K_n$, and
by Corollary 2.2,
 this means that
 $c_0(l_{p_n}^{K_n})_{n \in \N} \sim c_0(l_{q_n}^{K_n})_{n \in \N}.$ $\cqd$ 

\vs 0.3cm

{\bf Theorem 3.3:} {\it The relation $E_{K_{\sigma}}$ is Borel reducible to isomorphism and to complemented embeddability between subspaces of $c_0$.
Indeed there exist a sequence of integers $(K_n)_{n \in \N}$ and  a sequence of reals $(p_n)_{n \in \N}$ such that for $\alpha$ and $\beta$ in $X_0$, the following  statements are equivalent:}

(1) {\it $\alpha H_0 \beta.$} 

(2) {\it$c_0(l_{p_n}^{K_n}(\alpha)) \sim c_0(l_{p_n}^{K_n}(\beta)).$}

(3) {\it $c_0(l_{p_n}^{K_n}(\alpha)) \simeq c_0(l_{p_n}^{K_n}(\beta)).$}

(4) {\it $c_0(l_{p_n}^{K_n}(\alpha)) \st{c}{\hr} c_0(l_{p_n}^{K_n}(\beta)).$}

\vs 0.2cm

{\it Proof.} We choose $(K_n)_{n \in \N}$ and $(p_n)_{n \in \N}$ such that $p_1+1/\log K_1<2$, $p_n$ is decreasing to $1$, $n/\log K_n$ is decreasing, and
for all $n$, $p_n-p_{n+1} \geq 2n/{\log K_n}$.
This is possible if $\sum_{n=1}^{+\infty} \frac{n}{\log K_n}$
is small enough.
Then conditions (1), (2), (3) and (4) of Proposition 3.2 are achieved for any two sequences $(p_n+\frac{\alpha(n)}{\log K_n})_{n \in \N}$ and
$(p_n+\frac{\beta(n)}{\log K_n})_{n \in \N}$. 
Corollary 2.2 gives that (1) implies (2). Finally,  (4) implies (1) comes from Proposition 3.2. $\cqd$ 

\vs 0.3cm

{\bf Remark 3.4:} Observe that we cannot use Pe\l czy\'nski's decomposition method here to show
 that isomorphism and complemented bi-embeddability coincide, because the
 conditions we need to impose on $(p_n)_{n \in \N}$ and $(K_n)_{n \in \N}$ prevent the sequence
 $((p_n-p_{n+1}) \log K_n)_{n \in \N}$ from being bounded; that condition is needed to prove that   a $c_0$-sum of  $l_{p_n}^{K_n}$'s is isomorphic to its square.

We only used Banach spaces with unconditional bases. The
 crucial point in our method is that the spaces considered are isomorphic if and
only if their canonical bases are equivalent. As Rosendal proved that
 equivalence of Schauder bases is $K_{\sigma}$-complete \cite{R1}, we cannot hope to go further up in the hierarchy of complexity than
$K_\sigma$ with this method. So we now turn to a situation where isomorphism
corresponds to permutative equivalence of the canonical bases.

\vs 0.4cm

{\bf 4. Reducing $G$-equivalence rela\-tions to iso\-mor\-phism between
 separable Banach spaces}

\vs 0.2cm

 It is well-known that if
 $(Y_i)_{i \in \N}$ is  a sequence of Banach spaces,
 and  if a
Banach space $X$ is isomorphic to a subspace of $l_p(Y_i)_{i \in \N}$ for
some
$p \in [1,+\infty)$, then $X$ is isomorphic to a subspace of 
$\sum_{i=1}^n \oplus Y_i$ for some $n \in \N$ or $l_p$ is isomorphic to  a subspace of $X$.
In particular, if for some $p \in [1,+\infty)$, the space $l_p$ is isomorphic to a subspace of $l_{p_0}(l_{p_n})_{n \in \N}$, where $p_n \in [1,+\infty)$ for all $n \in \N \cup \{0\}$, then there exists $n \in \N \cup \{0\}$ such
 that  $p=p_n$. For a proof of these facts, see for example \cite{Bu} Theorem 1.1.

If $(X_n)_{n \in \N}$ is a sequence of Banach spaces, we shall define 
$l_p^{\infty}(X_n)_{n \in \N}$ as an $l_p$-sum where each $X_n$ appears in infinitely many summands.
In other words, $l_p^{\infty}(X_n)_{n \in \N} \simeq l_p(l_p(X_n))_{n \in \N}$,
where for each $n \in \N$, $l_p(X_n)$ denotes the $l_p$-sum of infinitely many copies of $X_n$, with permutative equivalence of the canonical bases.

\vs 0.3cm

{\bf Theorem 4.1:} {\it The relation $=^{+}$ is Borel reducible to isomorphism, to Lipschitz
isomorphism, to biembeddability  and to complemented biembeddability between
subspaces of $L_p$, $1 \leq p <2$.}

\vs 0.2cm

{\it Proof.} Fixing $1 \leq p <2$, we let $P$ be a perfect subset of the interval $]p,2[$.

We  define for $\alpha=(\alpha_n)_{n \in \N} \in P^{\omega}$, the Banach space
$$X(\alpha)=l_p^{\infty}(l_{\alpha_n})_{n \in \N}.$$
 \hs 0.6cm This defines a Borel map
and we show that it reduces $=^{+}$ on $P^{\omega}$ to isomorphism between
subspaces of $L_p$.
Indeed, first note that $X(\alpha)$ is isomorphic to a subspace of $L_p$ as an
$l_p$-sum of subspaces of $L_p$. 
Now if $\alpha =^+ \beta$, then every summand in the $l_p$-sum $X(\alpha)$
 (resp. $X(\beta)$) is a
summand in $X(\beta)$ (resp. $X(\alpha)$), and both appear infinitely many times as summands.
So $X(\alpha)$ is isometric to  
 $X(\beta)$ (in fact its canonical basis is permutatively equivalent to the
 canonical basis of $X(\beta)$).

Conversely, assume $X(\alpha)$ embeds in $X(\beta)$. Let $n \in \N$, we see
 that $l_{\alpha_n}$ is isomorphic to some subspace of  the $l_p$-sum $X(\beta)$. As $\alpha_n \neq p$,
 it follows that there exists $m$, such that $\alpha_n = \beta_m$ for some $m$.
Assuming $X(\beta)$ embeds in $X(\alpha)$, we get that $\beta_n = \alpha_{q}$ for some $q$.
As $n$ was arbitrary, we conclude that $\alpha =^+ \beta$.

Finally we conclude that for $\alpha$ and $\beta$ in $P^{\om}$,
$\alpha =^{+} \beta$ if and only if $X(\alpha)$ is isometric to $X(\beta)$, resp. isomorphic to, Lipschitz isomorphic to,
 complementably beimbeddable in, beimbeddable in $X(\beta)$. Once again we
 used \cite{HM} Theorem 2.4, together with reflexivity, to see that Lipschitz equivalence
 implies complemented biembeddability.  \pff $\cqd$ 

\vs 0.3cm

Before the next result we recall that  an operator $T$ from a Banach space
$X$ into a Banach space $Y$ is strictly singular if there exists no infinite
dimensional subspace $Z$ of $X$ such that the restriction of $T$ to $Z$ is an
isomorphism onto the image. Two Banach spaces $X$ and $Y$ are said to be totally incomparable if $X$ and $Y$ have no isomorphic closed subspaces of infinite dimension.

\vs 0.2cm

{\bf Theorem 4.2: (Wojtasczyk \cite{W})} {\it Assume that $X_1$ and $X_2$ are Banach spaces such that any operator from $X_1$ to $X_2$ is strictly singular. Let $X$ be a complemented subspace of $X_1 \oplus X_2$. Then $X$ is isomorphic to $Y_1 \oplus Y_2$, where $Y_i$ is a complemented subspace of $X_i$ for $i=1,2$.}

\vs 0.3cm

Given $R$ (resp. $R'$) an equivalence relation on a set $E$ (resp. $E'$), the
product $R \otimes R'$ is defined on $E \times E'$ by
$$(x,x')\ R \otimes R' \ (y,y') \Leftrightarrow x R x' \wedge y R y'.$$

{\bf Theorem 4.3:} {\it 
The relation
 $E_{K_{\sigma}} \otimes =^{+}$ is Borel reducible 
 to iso\-mor\-phism, Lips\-chitz iso\-mor\-phism and com\-ple\-men\-ted
 bi\-embed\-da\-bi\-lity bet\-ween subspaces of $L_p$, $1 \leq p <2$.}

\vs 0.1cm

{\it Proof.}
Fix $1 \leq p <2$. Let $f$ be a map given by Theorem 2.6 which Borel reduces
 $E_{K_{\sigma}}$
to isomorphism between subspaces of $l_{(p+1)/2}$. Let $P$ be a perfect subset
of $](p+1)/2,2[$  and $g$ be a map given by
Theorem 4.1 which Borel reduces
$=^+$ on $P^{\om}$
  to
  isomorphism between $l_p$-sums of
$l_{q_n}$-spaces for sequences $(q_n)$ in $P$.
 By the result at the beginning of this section, $l_{(p+1)/2}$ is totally incomparable with such $l_p$-sums.
Using Theorem 4.2, we check that the direct
sum $h$ of the two maps
 (defined by $h(\alpha,\beta)=f(\alpha) \oplus g(\beta)$) Borel reduces
 $E_{K_{\sigma}} \otimes =^{+}$ to isomorphism between
 subspaces of $L_p$.

Indeed, first note that by construction $h(\alpha,\beta)$ is a subspace of $L_p$, for
$\alpha$ in $X_0$ and $\beta$ in $P^{\om}$.
 Then assume $\alpha$ and $\alpha'$ in $X_0$ are $H_0$-related, and
$\beta$ and $\beta'$ in $P^{\om}$ satisfy $\beta =^+ \beta'$; then
 $h(\alpha,\beta) \simeq h(\alpha',\beta')$.
Conversely, assume $h(\alpha,\beta) \st{c} {\simeq} h(\alpha',\beta')$.
Then in particular, $g(\beta)$ is isomorphic to a complemented subspace of  $f(\alpha') \oplus
g(\beta')$. By  Theorem 4.2, it follows that $g(\beta) \simeq
U \oplus V$, with $U \st{c}{\hr} f(\alpha')$ and $V \st{c}{\hr} g(\beta')$. By
total incomparability of $g(\beta)$ and $f(\alpha')$, $U$ is finite dimensional. It follows that $g(\beta)
\st{c}{\hr} U \oplus g(\beta') \simeq g(\beta')$.
Symmetrically $g(\beta') \st{c}{\hr} g(\beta)$ and by Theorem 4.1, we deduce
that
$\beta=^+ \beta'$.

Similarly we get that $f(\alpha) \st{c}{\hr} f(\alpha') \oplus F$, where $F$
is finite-dimensional, as well as the complemented embedding in the other direction. We may then apply  Proposition
2.5 and finally get that $\alpha H_0 \alpha'$.

The claimed result is then obtained as before by  circular implications and \cite{HM}.
$\cqd$

\vs 0.4cm

{\bf 5. Final remarks, open problems and a conjecture.}

\vs 0.2cm

{\bf Remark 5.1:} Note that a Banach space not containing $l_2$ and without type $p$ for some $1
\leq p <2$ (resp. without cotype $q$ for some
$q>2$) has at least $3$ mutually non-isomorphic subspaces. Indeed: by
Gowers' dichotomy theorem \cite{Go1} 
and the fact that H.I. spaces are ergodic \cite{R1}, we may 
assume that there exists a subspace $X_1$ with an unconditional basis. By the
previously mentioned theorem of Komorowski and
Tomczak-Jaegermann \cite{KN}, some subspace $X_2$ of $X_1$ does not have an unconditional
basis, but has a basis (or at least a FDD in the case when $X$ does not have
non-trivial cotype), from which it follows that it has the approximation
property \cite{LT}. Finally the assumption about the type (resp. the cotype) and the results of
A. Szankowski \cite{S} imply the existence of a subspace $X_3$ without the
approximation property.

As a consequence of a study of subspaces of a Banach space with $k$-dimensional
 unconditional structure, for $k \in \N$, R. Anisca proved that for $X$ non isomorphic to $l_2$ and with
 finite cotype, $l_2(X)$ has countably mutually non-isomorphic subspaces \cite{A}.

 Finally it is proved in \cite{FR2} that every Banach space contains a subspace
which is a minimal space (that is, embeds in any of its subspaces) or
 contains continuum many mutually non isomorphic subspaces.

 \vs 0.2cm

With Remark 5.1, our results and those mentioned in 1.3 of the introduction,  known facts about complexity of isomorphism between
subspaces of a given Banach space may be seen in Figure 2.
For each equivalence relation $E$, we write the Banach spaces $X$ for which we
know that $E$ is Borel reducible to isomorphism between subspaces of $X$.

\

\begin{figure}[htbp]
\begin{center}
\psfrag{a}{{\scriptsize {$(1, =)$: $l_2$.}}}
\psfrag{b}{{\scriptsize {$(2, =)$: spaces not isomorphic to $l_2$.}}}
\psfrag{c}{{\scriptsize {$(3, =)$: spaces not containing $l_2$, without cotype $q$, for some}}}
\psfrag{r}{{\scriptsize { $q>2$ or without type $p$,  for some $1\leq p < 2$.}}}
\psfrag{n}{{\scriptsize {$(\omega, =)$: $l_2(X)$ for $X \not\simeq l_2$ with finite cotype.}}}
\psfrag{d}{{\scriptsize {$(2^{\omega}, =)$: spaces without a minimal subspace.}}}
\psfrag{e}{{\scriptsize {$E_{0}$: HI spaces, spaces with an unconditional basis}}}
\psfrag{E}{{\scriptsize {not isomorphic to their squares}}}
\psfrag{s}{{\scriptsize {(resp. to their hyperplanes),...}}}
\psfrag{f}{{\scriptsize {$E^{\infty}_{F_2}$}}}
\psfrag{g}{{\scriptsize {$=^{+}$}}}
\psfrag{h}{{\scriptsize {$E^{\infty}_{S_{\infty}}$}}}
\psfrag{i}{{\scriptsize {$E^{\infty}_{G_0}$}}}
\psfrag{j}{{\scriptsize {$E_{\Sigma^{1}_{1}}$}}}
\psfrag{k}{{\scriptsize {$E_{K_\sigma}$: $c_0$ and $l_p$, $1 \leq p <2$.}}}
\psfrag{l}{{\scriptsize {$E_{1}$: Tsirelson's space $T$.}}}
\psfrag{m}{{\scriptsize {$E_{K_\sigma} \otimes=^{+}$: $L_p$, $1 \leq p <2$.}}}
\psfrag{q}{{\scriptsize {\hbox{Figure 2: diagram  of complexity of isomorphism
        between subspaces of a separable Banach space.}}}}
\includegraphics[scale=0.35]{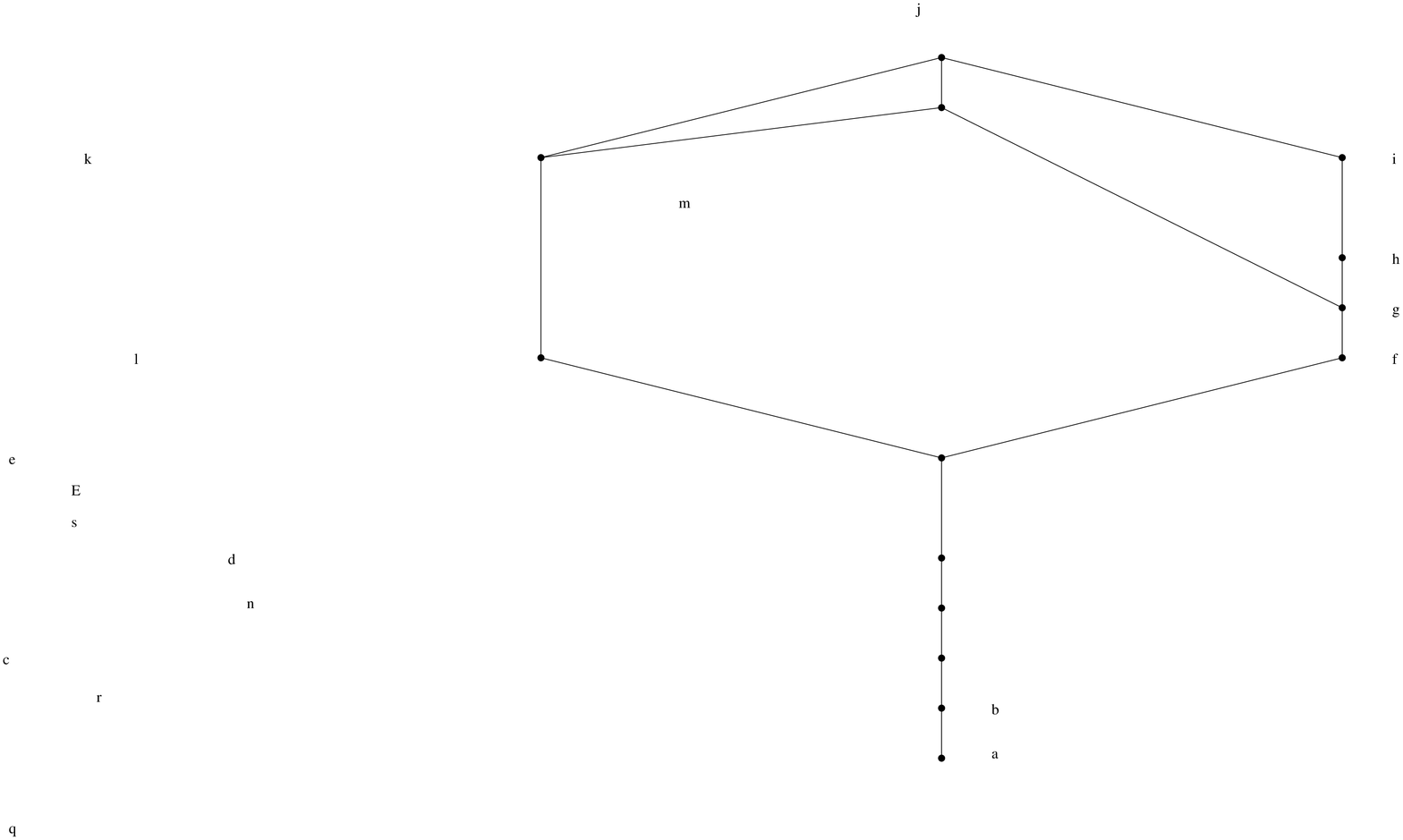}
\label{figura2006}
\end{center}
\end{figure}

{\bf Remark 5.2:} A problem we left open is whether we may extend our results to prove
that spaces
$l_p$ are ergodic for $p>2$. Our results about $c_0$ and
$l_p, 1 \leq p <2$
suggest the conjecture that $l_2$ is the only non-ergodic Banach
space. It is also of interest to restrict the question to
particular subspaces, such as block-subspaces of a given basis. As any
normalized block-basis of the canonical basis of $c_0$ or $l_p$ is equivalent to the
original basis, these spaces would be, as is maybe natural, of the lowest
complexity possible. We conjecture the following. 

\vs 0.2cm

{\bf Conjecture 5.3:}  Let $X$ be a Banach space with  an unconditional basis.
  Then either the relation $E_0$ is Borel reducible to isomorphism between
 block-sub\-spaces of $X$, or $X$ is isomorphic to
 $c_0$ or $l_p$, $1 \leq p <+\infty$.

\vs 0.2cm

{\bf Remark 5.4:} The question of the exact complexity of isomorphism between
 separable Banach spaces (seen as subspaces of $C([0,1])$ with the
 Effros-Borel structure), or even between Banach spaces with  Schauder bases
 (seen as subsequences of the universal basis of Pe\l czy\'nski) is quite open. It could be
 $\Sigma_1^{1}$-complete, that is, $\leq_B$-maximum among analytic equivalence
 relations.

The limitation for our methods might come from the fact that our examples are
isomorphic exactly when their canonical bases are permutatively equivalent. If
the complexity of permutative equivalence
 of basic sequences is too low, it will
be necessary to find quite different types of reductions. It could be
interesting to work with general Banach lattices instead of discrete ones.

\vs 0.2cm

{\bf Remark 5.5:} Gowers \cite{Go} solved the so-called Schroeder-Bernstein
 problem for Banach spaces, by proving  that complemented biembeddability and
 isomorphism between Banach spaces need not coincide in general. We notice  that for our
 examples, they do coincide. This is probably because many techniques known about isomorphic spaces concern, more generally, complemented subspaces. We may wonder how far these two properties are from each other from a point of view of complexity.
In this direction, we show in \cite {FG} how to construct a continuum of
 mutually non-isomorphic subspaces which are however complemented in each other.

\vs 0.2cm
{\bf Acknowledgments.} We thank C. Rosendal for answering numerous 
questions and doubts about the theory of Borel reducibility between analytic
equivalence relations on Polish spaces. We also thank G. Godefroy and N. Kalton 
for drawing our attention to the results of Casazza-Kalton in \cite{CK1} and  \cite{CK2}.

\vs 0.3cm

 Equipe d'Analyse Fonctionnelle,

Universit\'e Paris 6,

Bo\^ite 186, 4, Place Jussieu,

 75252, Paris Cedex 05,

 France.

 E-mail: ferenczi@ccr.jussieu.fr.

\

and

\

Department of Mathematics, IME.

University of S\~ao Paulo.

05311-970 S\~ao Paulo, SP,

Brasil.

E-mail: ferenczi@ime.usp.br, eloi@ime.usp.br

\end{document}